\documentclass{amsart}
\usepackage{amssymb}
\usepackage{amsthm}
\usepackage{amscd}
\usepackage{graphics}
\newtheorem*{reftheorem}{Theorem}

\newtheorem*{theorem 1}{Theorem 1}
\newtheorem*{theorem 2}{Theorem 2}
\newtheorem*{theorem 3}{Theorem 3}
\newtheorem{lemma}{Lemma}[section]

\newtheorem{prop}{Proposition}[section]
\newtheorem{introtheorem}{Theorem}
\newtheorem{introlemma}{Lemma}
\newtheorem{introprop}{Proposition}

\newtheorem*{refprop}{Proposition}

\newtheorem*{conj}{Conjecture}

\newtheorem*{iso}{Isoperimetric Inequality in an Ellipse}
\newtheorem*{diriso}{Directional Isoperimetric Inequality in an Ellipse}
\numberwithin{equation}{section}
\title{Directional isoperimetric inequalities and rational homotopy invariants}
\author{Larry Guth}
\address{Department of Mathematics, Stanford, Stanford CA, 94305 USA}
\email{lguth@math.stanford.edu}

\begin{document}
\begin{abstract} We estimate the second order linking invariants
of Lipschitz maps from an n-dimensional ellipse.  The estimate
uses a new directionally-dependent version of the isoperimetric
inequality for cycles inside the ellipse.  Using this work,
we prove new lower bounds for the k-dilation of maps from one
ellipse to another.
\end{abstract}

\maketitle

In this paper, we estimate a second-order rational homotopy
invariant of a map in terms of the map's Lipschitz constant. 
This problem turns out to be qualitatively harder than estimating
a first-order rational homotopy invariant such as the Hopf
invariant.

In \cite{Gr1} and \cite{Gr2}, Gromov described a basic upper
bound for the rational homotopy invariants of a Lipschitz map
from a Riemannian manifold.  In \cite{G1}, I showed that when the
domain is an n-dimensional ellipse, then the estimates for the
Hopf invariant and the linking invariant are sharp up to a
constant factor.  (Recall that an n-dimensional ellipse $E$ with
principal axes $E_0 \le ... \le E_n$ is the set $\{ x \in
\mathbb{R}^{n+1} | \sum_{i=0}^n (x_i / E_i)^2 = 1 \}$.)

We study a second-order linking invariant of maps from $S^n$ to a
wedge of three spheres $S^{k_1} \vee S^{k_2} \vee S^{k_3}$. 
Gromov's method gives an upper bound for this invariant in terms
of a metric on the domain, a metric on the range, and the
Lipschitz constant of the map.  This upper bound, however, may be
too large.  As we will see, even if the domain is an ellipse
and each sphere in the range has the unit sphere metric, the
upper bound may be much too large.  In this special case, we will
give a better upper bound which is sharp up to a constant factor.

Before I state the results, I want to say something about the new
method.  Gromov's method uses the isoperimetric inequality.  The
new idea in this paper is to replace ordinary isoperimetric
inequalities by directional isoperimetric inequalities that
separately keep track of the amount of volume of a surface
pointing in different directions.  If $J$ is an m-tuple of
integers from 1 to n, let $P(J)$ denote the corresponding
coordinate m-plane in $\mathbb{R}^n$.  Now, if $C$ is an
m-dimensional surface in $\mathbb{R}^n$, then we define
$Vol_J(C)$ to be the volume of the projection of $C$ onto the
corresponding m-plane.  For example, suppose that $C$ is the long
thin curve in the figure below.

\includegraphics{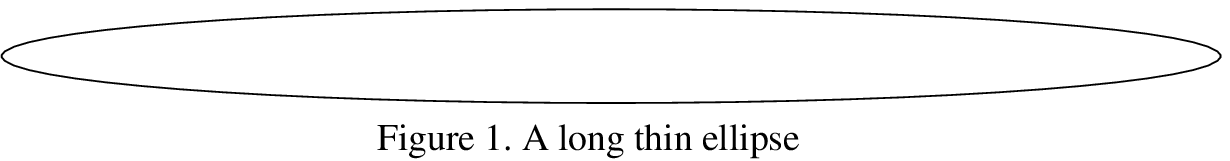}

\noindent Suppose the curve $C$ has 
$Vol_1(C) = 20$ and $Vol_2(C) = 4$.  The total
length of $C$ is slightly more than 20.  In
general, a plane curve of length 20 may bound a region
of area 30, but because this curve has such small
2-volume, it bounds a significantly smaller region.  By keeping
track of the volumes in different directions, we will be able to
find more efficient fillings of certain cycles.

Loomis and Whitney proved a directional estimate in this spirit
in \cite{LW}.

\begin{reftheorem} (Loomis, Whitney) Suppose that $U \subset
\mathbb{R}^n$ is an open set.  For any (n-1)-tuple $J$, let
$A_J$ denote the area of the projection of $U$ onto $P(J)$.
Then the volume of $U$ is bounded in terms of the areas of
its projections as follows.

$$Vol(U) \le [\prod_J A_J]^{\frac{1}{n-1}}.$$

In this formula, $J$ varies over the $n$ different (n-1)-tuples
of integers from 1 to n.
\end{reftheorem}

The Loomis and Whitney theorem bounds the volume enclosed by a hypersurface
in terms of the directional volumes of the hypersurface.  In this paper
we derive estimates in a similar spirit which apply to surfaces of
any codimension.

Now we return to estimating rational homotopy invariants of
Lipschitz maps.  In the rest of the introduction, I will try to
explain the answers to three questions.

1. How can we estimate rational homotopy invariants using
isoperimetric inequalities?

2. Why are these estimates far from sharp for some second-order
invariants?

3. How can we improve the estimates using directional
isoperimetric inequalities?

The linking invariant is defined for a map $F$ from $S^n$ to the
wedge of spheres $S^{k_1} \vee S^{k_2}$ provided that the
dimensions obey $n + 1 = k_1 + k_2$ and $2 \le k_1 \le k_2$.  We
let $q_1$ denote a generic point in $S^{k_1}$ and $z_1$ denote
the fiber $F^{-1}(q_1)$.  The fiber $z_1$ is a closed oriented
submanifold of $S^n$.  We let $y_1$ be a chain in $S^n$ with
$\partial y_1 = z_1$.  We
let $q_2$ be a generic point of $S^{k_2}$, and we define $z_2$ to
be the intersection $F^{-1}(q_2) \cap y_1$.  For generic
$q_2$, this intersection
will consist of finitely many points, each with an orientation. 
The signed number of points is the linking invariant of $F$,
denoted $L(F)$.  (By standard results in topology, $L(F)$ does
not depend on the choices we made, and it is a homotopy
invariant of $F$.)

Now suppose that we have a metric $g$ on $S^n$ and also a metric
on the target $S^{k_1} \vee S^{k_2}$.  Using these metrics, we
can define the Lipschitz constant of the map $F$.  Our goal is to
bound the linking invariant of $F$ in terms of the Lipschitz
constant and the metrics.  In particular, I want to understand
how the best bound depends on the metric $g$. In order to control
the linking invariant of a map $F$, we will estimate the number
of points in $z_2$, and in order to do that, we control
the volumes of $z_1$ and $y_1$.  Our estimate involves two
geometric ingredients.  The first ingredient is to find a fiber
with small volume, using the coarea inequality.

\begin{introlemma} (Coarea inequality) 
Let $F: (M^n, g) \rightarrow (N^q, h)$ be a
$C^\infty$ map with Lipschitz constant $L$.  Then $F$ has a
regular fiber with volume at most $L^q Vol(M) / Vol(N)$.
\end{introlemma}

Using this lemma, we can bound the volume of $z_1$ in terms of
the volume of $(S^n, g)$.  The next step of the argument is to
bound the volume of $y_1$ in terms of the volume of $z_1$.  This
step requires an isoperimetric inequality that holds in $(S^n,
g)$.  In this paper, we will focus on the case that $(S^n, g)$ is an
ellipse.  In that case, the relevant isoperimetric
inequality is described by the following lemma.

\begin{introlemma} (\cite{G1}) Suppose that $E$ is an 
n-dimensional ellipse with principal axes $E_0 \le ... \le E_n$. 
Suppose that $z$ is an m-cycle in $E$.  Then there is an
(m+1)-chain $y$ in $E$ with $\partial y = z$ obeying the following
estimate.

$$Vol(y) \le C(n) [E_{m+1} + E_{n-m}] Vol(z).$$

\end{introlemma}

Suppose that $F$ is a map from the ellipse $E$ to the
wedge $S^{k_1} \vee S^{k_2}$ equipped with its standard metric. 
We can bound the linking invariant of $F$ in terms of its
Lipschitz constant by combining the two lemmas above.  First, we
use Lemma 1 to bound the volume of $z_1$.  Then we use Lemma 2 to
bound the volume of $y_1$.  Finally, we use Lemma 1 again to
bound the volume of $z_2$.  Since $z_2$ is a 0-cycle, its volume
is equal to the number of points in it, and so the volume of
$z_2$ controls the linking invariant.  Putting together these
steps, we get the following estimate.

\begin{introprop} (\cite{G1}) Let $E$ be an n-dimensional ellipse
with principal axes $E_0 \le ... \le E_n$.  Let $2 \le k_1
\le k_2$ with $n+1 = k_1 + k_2$.  Let $F$ be a map from $E$ 
to $S^{k_1} \vee S^{k_2}$.  Equip the target with its
standard metric.  Suppose that $F$ has Lipschitz constant $L$.

$$\textrm{Then } |L(F)| \le C(n) E_{k_2} Vol(E) L^{n+1}.$$

\end{introprop}

\noindent In \cite{G1}, I showed that this estimate is sharp up
to a factor $C(n)$ for sufficiently large $L$.

Now we turn to second order invariants.  We will define a
second-order linking invariant.  It is closely analogous to
the linking invariant we just considered, but it involves three
fibers and two filling operations.

The second-order linking invariant is defined for a map $F$ from
$S^n$ to the wedge of spheres $S^{k_1} \vee S^{k_2} \vee S^{k_3}$
provided that the dimensions obey $n + 2 = k_1 + k_2 + k_3$ and
$2 \le k_1 \le k_2 \le k_3$.  We let $q_1$ denote a generic point
in $S^{k_1}$ and $z_1$ denote the fiber $F^{-1}(q_1)$.  We let
$y_1$ be a chain in $S^n$ with $\partial y_1 = z_1$.  We let
$q_2$ be a generic point of $S^{k_2}$, and we define $z_2$ to be
the intersection $F^{-1}(q_2) \cap y_1$.  The fiber $z_2$ is an
integral cycle in $S^n$.  We let $y_2$ be a chain in $S^n$ with
$\partial y_2 = z_2$.  Finally, we let $q_3$ be a generic point
of $S^{k_3}$, and we define $z_3$ to be the intersection
$F^{-1}(q_3) \cap y_2$.  This intersection will consist of
finitely many points, each with an orientation.  The signed
number of points is the second-order linking invariant of $F$,
denoted $L_2(F)$.  (Like the first-order linking invariant,
$L_2(F)$ does not depend on the choices we made, and it is a
homotopy invariant of $F$.)

We can use the same strategy to bound $L_2(F)$ for a Lipschitz
map from an ellipse.  Lemma 1 bounds the volume of $z_1$.  Then
Lemma 2 bounds the volume of $y_1$.  Then Lemma 1 bounds the
volume of $z_2$.  Then Lemma 2 bounds the volume of $y_2$. 
Finally, Lemma 1 bounds the volume of $z_3$ which bounds
$|L_2(F)|$.  If we carry out the calculations, we get
the following bounds.

\begin{introprop} Suppose that $2 \le k_1 \le k_2 \le k_3$ and
that $n + 2 = k_1 + k_2 + k_3$.  Let $E$ be an n-dimensional
ellipse with principal axes $E_0 \le ...
\le E_n$.  Let $F$ be a map from $E$ to the wedge of unit spheres
$S^{k_1} \vee S^{k_2} \vee S^{k_3}$ with Lipschitz constant
$L$.

If $k_3 < (n+1)/2$, then $L_2(F)$ is bounded as follows.

$$|L_2(F)| \le C(n) E_{n-k_1 + 1} E_{n-k_3 + 1} Vol(E) L^{n+2}.
$$

If $k_3 \ge (n+1)/2$, then $L_2(F)$ is bounded as follows.

$$|L_2(F)| \le C(n) E_{n-k_1 + 1} E_{k_3} Vol(E) L^{n+2}.
$$

\end{introprop}

In the first case, $k_3 < (n+1)/2$, it turns out that this
inequality is not sharp.  The main result of this paper is a
refined inequality which is sharp up to a constant factor.

In order to understand why this basic upper bound is not sharp, we
have to understand a little bit about the isoperimetric
inequality in ellipses given in Lemma 2 above.  According to Lemma 2,
an m-cycle $z$ bounds a chain $y$ with $|y| \lesssim [E_{m+1} +
E_{n-m}] |z|$.  This bound is sharp up to a constant factor.  The
way that the worst case cycle looks depends on the dimension $m$. 
If $m \ge (n-1)/2$, then the smallest m-dimensional equator of
$E$ is the hardest m-cycle to fill (up to a constant factor).  On
the other hand, if $m < (n-1)/2$, then the largest m-dimensional
equator of $E$ is the hardest m-cycle to fill.

If $k_3 < (n+1)/2$, then the cycle $z_1$ has dimension $m_1$ at
least $(n-1)/2$, but the cycle $z_2$ has dimension $m_2$ less
than $(n-1)/2$.  Now we can informally describe why Proposition 2
is not sharp for $k_3 < (n+1)/2$.  If Proposition 2 were sharp,
it would mean that $z_1$ ``looks like'' the
smallest $m_1$-dimensional equator of $E$ and $z_2$ ``looks
like'' the largest $m_2$-dimensional equator of $E$.  If $z_1$
actually were the smallest $m_1$-dimensional equator of $E$,
then we could choose $y_1$ to be a hemisphere inside the smallest
$m_1+1$-dimensional equator of $E$.  In that case, $z_2$ would
be an $m_2$-cycle lying inside the smallest $m_1+1$-dimensional
equator of $E$.  Such a cycle looks very different from 
the largest $m_2$-dimensional equator of $E$, and in particular
it can be filled much more efficiently.  We can make this
argument effective by keeping track of the directional volumes
of $y_i$ and $z_i$.

Now we describe the argument in a bit more detail.  We begin in
the same way as the basic argument: we apply Lemma 1 to find a
fiber $z_1$ with controlled volume.  In the basic argument, we
applied Lemma 2 to find a chain $y_1$ with volume at most $\sim
E_{m_1+1} |z_1|$.  In our more refined argument, we construct a
chain $y_1$ by a different method, which allows us to bound the
directional volumes of $y_1$ in a useful way.  In some
directions, the volume of $y_1$ may be as large as $E_{m_1+1}
|z_1|$, but in most directions the directional volume of $y_1$ is
much smaller.  The directional volume of $y_1$ is concentrated in
the directions where the ellipse $E$ is small, such as $I = [1,
..., m_1+1]$.  Next we choose a fiber $z_2 \subset y_1$.  In the
basic argument, we applied Lemma 1 to find a fiber $z_2$ with
controlled volume.  In the refined argument, we want to control
all the directional volumes of $z_2$ in terms of our bounds for
the directional volumes of $y_1$.  To do that, we use a small
generalization of Lemma 1.  At this stage, we have a bound for
the total volume of $z_2$ which is the same as in the basic
argument, but we have stronger bounds on many of the directional
volumes of $z_2$ which show that $z_2$ is concentrated in the
directions where $E$ is small.  Next, we again use a directional
isoperimetric inequality to find a chain $y_2$ with boundary
$z_2$.  The basic isoperimetric inequality in Lemma 2 tells us
that we can find $y_2$ with volume at most $\sim E_{n-m_2}
|z_2|$.  But since the volume of $z_2$ is concentrated in the
small directions, we can improve on that estimate and find $y_2$
with significantly smaller volume.  Finally we use Lemma 1 to
bound the volume of $z_3$ and thus bound $|L_2(F)|$ as before.

\begin{introtheorem} Suppose that $2 \le k_1 \le k_2 \le k_3$,
and $n + 2 = k_1 + k_2 + k_3$.  Let $E$ be an n-dimensional
ellipse with principal axes $E_0 \le ... \le E_n$.  Let
$F$ be a map from $E$ to the wedge of unit spheres $S^{k_1} \vee
S^{k_2} \vee S^{k_3}$ with Lipschitz constant at most $L$.  Then
$L_2(F)$ is bounded as follows.

$$|L_2(F)| \le C(n) E_{n-k_1 + 1} E_{k_3} Vol(E) L^{n+2}.
$$

\end{introtheorem}

If $k_3 < (n+1)/2$, then the upper-bound in Theorem 1 is
better than the one in Proposition 2 by a factor $E_{n-k_3+1} /
E_{k_3}$, which may be arbitrarily large.  On the other hand, the
upper bound in Theorem 1 is sharp up to a constant factor.

\begin{introtheorem} Suppose that $2 \le k_1 \le k_2 \le k_3$,
and $n + 2 = k_1 + k_2 + k_3$.  Let $E$ be an n-dimensional
ellipse with principal axes $E_0 \le ...
\le E_n$.  Suppose that $L > C(n) E_1^{-1}$.  Then there is a map
$F$ from $E$ to the wedge of unit spheres
$S^{k_1} \vee S^{k_2} \vee S^{k_3}$ with Lipschitz constant $L$
and $L_2(F)$ bounded below as follows.

$$L_2(F) \ge c(n) E_{n-k_1 + 1} E_{k_3} Vol(E) L^{n+2}.
$$

\end{introtheorem}

We pause here to make some comments.  For a general Riemannian
metric $(S^n, g)$, directional volumes are not even well-defined. 
The refined estimate in Theorem 1 is applicable only to
ellipsoidal metrics, whereas the basic estimate is applicable to
all metrics.  But in its narrow range of applicability, the
refined estimate outperforms the basic estimate and is sharp up
to a constant factor.

I became interested in this question because I was trying to
estimate the k-dilations of mappings from one ellipse to another.
Recall that the k-dilation is a generalization of the Lipschitz
constant that measures how much a mapping stretches k-dimensional
areas.  Gromov noticed that his upper bounds for maps with a given
Lipschitz constant extend to maps with a given k-dilation for
an appropriate value of $k > 1$ depending on the problem.  Similarly,
Theorem 1 can be extended to maps with a bound on the k-dilation.
More precisely, if  $F$ has $k_1$-dilation at most $L^{k_1}$,
then the conclusion of Theorem 1 still holds.  Hence, we can
estimate the largest value $L_2(F)$ for a map $F$ from $E$ to
a standard wedge of spheres with a given $k_1$-dilation.  This
result implies some new estimates about the k-dilations of maps
from one ellipse to another.

\begin{introtheorem} Let $E, E'$ be n-dimensional ellipses. 
Let $E_0 \le ... \le E_n$ be the principal axes of $E$.  Let
$E_0' \le ... \le E_n'$ be the principal axes of $E'$.  Let $Q_i
= E_i'/E_i$.  Suppose that $\Phi$ is a map from $E$ to $E'$ with
degree $D$.  Suppose that $2 \le k_1 \le k_2 \le k_3$, $n+2 = k_1
+ k_2 + k_3$ and $k \le k_1$.  Then the following inequality
holds.

$$Dil_k(\Phi) > c(n) [|D| Q_{n-k_1+1} Q_{k_3} Q_1 ...
Q_n]^{\frac{k}{n+2}}.$$

\end{introtheorem}

The simplest example is the case $n=4$ and $k=k_1=k_2=k_3=2$.  In
this case, the 2-dilation of $\Phi$ is at least $\sim |D|^{1/3}
Q_1^{1/3} Q_2^{2/3} Q_3^{2/3} Q_4^{1/3}$.

The problem of giving sharp lower bounds for the k-dilation of a
degree 1 map from one n-dimensional ellipse to another looks very
difficult.  (To be clear, I would like an estimate which is sharp
up to a constant factor $C(n)$ independent of the principal axes
of the ellipses.)  The analogous problem for the 2-dilation of
mappings between 4-dimensional rectangles was recently solved in
\cite{G4}, and the answer turned out to be complicated.  The
near-optimal mappings are far from linear.  The possible pairs of
rectangles are divided into several cases and in each case there
is a rather different non-linear mapping.  Also, it turns out
that the smallest 2-dilation of a degree 1 diffeomorphism may be
larger than the smallest 2-dilation of a degree 1 map by an
arbitrary factor.

\vskip5pt

{\it Acknowledgements.} This paper is based on a section of my
thesis.  I am grateful to my thesis advisor, Tom Mrowka, for his
help and support.

\section{Background}

We recall the definition of $L_2(F)$ and prove that it is a homotopy
invariant.  The invariant $L_2(F)$ is a special case of a rational
homotopy invariant and so it fits into a general theory of rational
homotopy invariants and differential forms developed by Sullivan
\cite{S}.  (Historical question: who first defined the invariant $L_2$?)
I don't know a good reference in the literature, so we give a
self-contained presentation along the lines of Bott and Tu
(\cite{BT}).  We define the invariant $L_2(F)$ and prove that it
is a homotopy invariant.

Suppose that $F$ is a map from $S^n$ to the wedge of spheres $S^{k_1} \vee
S^{k_2} \vee S^{k_3}$.  We want to define $L_2(F)$ in the case
that the dimensions obey the conditions $2 \le k_1 \le k_2 \le
k_3$ and $n+2 = k_1 + k_2 + k_3$.  

Let $q_i$ denote a point in $S^{k_i}$.  First we consider the
special case that $F$ is a $C^\infty$ map and that $q_i$ are
regular values for $F$.  (We mean here that $F$ is $C^\infty$
away from the inverse image of the basepoint of $S^{k_i}$.)  In
this case, we define $L_2(F)$ as follows.  First we consider the
fiber $z_1 = F^{-1}(q_1)$ in $S^n$.  This fiber is an orientable
manifold of dimension $n - k_1$.  We choose an integral chain
$y_1 \subset S^n$ with $\partial y_1 = z_1$.  Next we consider
the intersection of $y_1$ with the fiber $F^{-1}(q_2)$.  After
putting $y_1$ in general position, the intersection is a cycle
$z_2$ of dimension $n - k_1 + 1 - k_2 = k_3 - 1$.  We choose a
chain $y_2$ in $S^n$ with $\partial y_2 = z_2$.  The chain $y_2$
has dimension $n - k_1 - k_2 + 2 = k_3$.  We let $\pi$ denote the
retraction from $S^{k_1} \vee S^{k_2} \vee S^{k_3}$ to $S^{k_3}$. 
The compositon $\pi \circ F$ maps $y_2$ to $S^{k_3}$ and maps the
boundary of $y_2$ to the basepoint of $S^{k_3}$.  Therefore, the
composition has a well-defined degree, and the degree is defined
to be $L_2(F)$.

For a general map $F$, we homotope $F$ to a map $F_{nice}$ which
is $C^\infty$ and for which $q_i$ are regular values.  Then we
define $L_2(F)$ to be $L_2(F_{nice})$.

We have to prove that $L_2(F)$ is well-defined, and does not
depend on the choices that we made above.  In particular, it does
not depend on the choice of $y_i$, it does not depend on the
choice of $q_i$, and it does not depend on the choice of
$F_{nice}$.  Also, we prove that $L_2$ is a homotopy invariant of $F$.

\begin{prop} The quantity $L_2(F)$ is independent of the choices.
First, if $F$ is $C^\infty$ and $q_i$ are regular values of $F$,
then $L_2(F)$ is independent of the choices of $y_1$ and $y_2$. 
Second, $L_2(F)$ is independent of the choice of $F_{nice}$. 
Third, $L_2$ is a homotopy invariant.  Finally, $L_2(F)$ is
independent of the choice of $q_i$.
\end{prop}

\proof We assume that $F$ is $C^\infty$ and that $q_i$ are
regular values of $F$.  A priori, $L_2$ depends on $F$, $y_1$,
and $y_2$, so we write it as $L_2(F, y_1, y_2)$.  We prove that
$L_2(F)$ is independent of $y_1$ and $y_2$.

First we show that $L_2$ is independent of the choice of $y_2$.
Suppose we chose a different cycle $y_2'$.  Let $\Sigma$ denote
$y_2 - y_2'$.  Note that $\Sigma$ is a $k_3$-cycle in $S^n$.  The
difference $L_2(F, y_1, y_2) - L_2(F, y_1, y_2')$ is given by
the degree of the cycle $\pi \circ F (\Sigma)$ in $S^{k_3}$.  But
the cycle $\Sigma$ is exact, so this degree is zero.  Since $L_2$
does not depend on the choice of $y_2$, we may write it as
$L_2(F, y_1)$.

Next we show that $L_2$ is independent of the choice of $y_1$.  We
study $L_2(F, y_1) - L_2(F, y_1')$.  The difference $y_1 - y_1'$ is
a cycle in $S^n$ that bounds a chain $A$.  Now we define $z_2$
to be the intersection of $y_1$ with $F^{-1}(q_2)$ and $z_2'$ to
be the intersection of $y_1'$ with $F^{-1}(q_2)$.  We define $B$
to be the intersection of $A$ with $F^{-1}(q_2)$.  The boundary
of $B$ is $z_2 - z_2'$.  Next we pick a chain $y_2$ with boundary
$z_2$, and we define $y_2'$ to be $y_2 - B$.  Note that the
boundary of $y_2'$ is $z_2 - z_2 + z_2' = z_2'$.  Recall that
$L_2(F, y_1)$ is the degree of $\pi \circ F(y_2)$ and $L_2(F, y_1')$
is the degree of $\pi \circ F(y_2')$.  Therefore the difference
$L_2(F, y_1) - L_2(F, y_1')$ is the degree of $\pi \circ F(B)$.  But
$B$ lies in the fiber $F^{-1}(q_2)$ and $q_2$ is a point in $S^{k_2}$,
so $\pi \circ F(B)$ is the basepoint $*$.  Since $L_2$ does not
depend on the choice of $y_1$, we may write it as $L_2(F)$.

Now we suppose that $F$ is a homotopy from $F_0$ to $F_1$.  We
assume that $q_i$ are regular values for $F$, $F_0$, and $F_1$. 
Under these hypotheses, we prove that $L_2(F_0) = L_2(F_1)$.

Suppose that $F: S^n \times [0, 1] \rightarrow S^{k_1} \vee S^{k_2}
\vee S^{k_3}$ is a homotopy from $F_0$ to $F_1$.  We have to check
that $L_2(F_0) = L_2(F_1)$.  First we consider the fiber
$F^{-1}(q_1)$. This fiber is a homology from $F_0^{-1}(q_1)$ to
$F_1^{-1}(q_1)$.  We let $y_{1, 0}$ be a chain filling
$F_0^{-1}(q_1)$ in $S^n \times \{ 0 \}$, and we let $y_{1, 1}$ be
a chain filling $F_1^{-1}(q_1)$ in $S^n \times \{ 1 \}$.  The sum
$y_{1,1} + F^{-1}(q_1) - y_{1,0}$ defines a cycle in $S^n \times
[0,1]$, and we defined $y_1$ to be a chain filling this cycle. 
The dimension of $F^{-1}(q_1)$ is $n+1 - k_1$, and so the
dimension of $y_1$ is $n+2 - k_1$.

Next we intersect $y_1$ with the fiber $F^{-1}(q_2)$.  We make
the following definitions.

$z_{2, 0} := y_{1,0} \cap F^{-1}(q_2) \subset S^n \times \{ 0 \}.$

$z_{2, 1} := y_{1,1} \cap F^{-1}(q_2) \subset S^n \times \{ 1 \}.$

$z_2 := y_1 \cap F^{-1}(q_2) \subset S^n \times [0, 1].$

Here $z_{2, 0}$ and $z_{2,1}$ are cycles, and $z_2$ is a chain
with boundary $z_{2,1} - z_{2,0}$.

Next we fill $z_2$.  We let $y_{2,0}$ be a chain in $S^n \times
\{ 0 \}$ with boundary $z_{2,0}$.  We let $y_{2,1}$ be a chain in
$S^n \times \{ 1 \}$ with boundary $z_{2,1}$.  And we let $y_2$
be a chain in $S^n \times [0,1]$ with boundary $y_{2,1} + z_2 -
y_{2,0}$.

Finally, we consider the map $\pi \circ F$ from $y_2$ to
$S^{k_3}$.  The map $\pi \circ F$ takes $z_2$ to the basepoint. 
The image $\pi \circ F (y_{2,0})$ is a $k_3$-cycle in $S^{k_3}$
of degree $L_2(F_0)$.  The image $\pi \circ F(y_{2,1})$ is a
$k_3$-cycle in $S^{k_3}$ of degree $L_2(F_1)$.  The image $\pi
\circ F (y_2)$ is a homology from $\pi \circ F (y_{2,0})$ to $\pi
\circ F (y_{2,1})$.  Hence $L_2(F_0) = L_2(F_1)$. 

Now for any map $F$, we can homotope $F$ to a $C^\infty$ map
$F_{nice}$ for which $q_i$ are regular values.  We define
$L_2(F)$ to be $L_2(F_{nice})$.  Because of the homotopy result
we proved above, the value of $L_2(F)$ does not depend on how we
choose $F_{nice}$.  It follows that $L_2$ is a homotopy invariant
of $F$.  

Finally, we check that $L_2$ does not depend on the choice of
$q_i$ as long as $q_i \in S^{k_i}$ and $q_i$ is not the base
point. Let $\tilde q_i \in S^{k_i}$ be some other points, and let
$\tilde L_2$ be the linking invariant defined using $\tilde q_i$
in place of $q_i$. Let $G$ be a diffeomorphism of $S^{k_1} \vee
S^{k_2} \vee S^{k_3}$, homotopic to the identity, taking $\tilde
q_i$ to $q_i$.  Let $F$ be any map from $S^n$ to $S^{k_1}
\vee S^{k_2} \vee S^{k_3}$ so that $\tilde q_i$ are regular
values.  Then $G \circ F$ has $q_i$ as regular values, and
$L_2(G \circ F) = \tilde L_2(F)$.  But since $L_2$ is a homotopy invariant
and $G$ is homotopic to the identity, $L_2(G \circ F) = L_2(F)$. \endproof

Next we consider an example to show that the $L_2$ invariant can
be non-trivial.  Suppose that $f: S^n \rightarrow S^{k_1} \vee
S^{n-k_1+1}$ is a continuous map.  Suppose that $g: S^{n-k_1+1}
\rightarrow S^{k_2} \vee S^{k_3}$ is a continuous map.  We assume
as usual that $n+2 = k_1 + k_2 + k_3$ and $2 \le k_1 \le k_2 \le
k_3$, and therefore the linking invariants of $f$ and $g$ are
each defined.  We let $g^+$ denote the map from $S^{k_1} \vee
S^{n-k_1+1}$ to $S^{k_1} \vee S^{k_2} \vee S^{k_3}$ which is
equal to the identity on $S^{k_1}$ and is equal to $g$ on
$S^{n-k_1+1}$. Then the composition $g^+ \circ f$ maps $S^n$ to
$S^{k_1} \vee S^{k_2} \vee S^{k_3}$.  The second-order linking
invariant $L_2(g^+ \circ f)$ is equal to the product $L(g) L(f)$,
which may be non-zero.

For completeness, we calculate $L_2(g^+ \circ f)$.  Let $F = g^+
\circ f$.  We let $z_1 = F^{-1}(q_1) = f^{-1}(q_1)$.  Then we
choose a chain $y_1$ with $\partial y_1 = z_1$.  Next we let $z_2
= y_1 \cap F^{-1}(q_2) = y_1 \cap f^{-1}[g^{-1}(q_2)]$. Now we
let $w$ be $g^{-1}(q_2)$, which is a cycle in $S^{n-k_1+1}$.  We
choose a chain $v$ with $\partial v = w$.  Next we have to choose
a chain $y_2$ with $\partial y_2 = z_2$.  The trick in this
calculation is that we choose $y_2 = y_1 \cap f^{-1}(v)$. 
Finally, we define $z_3 = y_2 \cap F^{-1}(q_3)$.  Expanding this
formula, we see $z_3 = y_1 \cap f^{-1}(v) \cap
f^{-1}[g^{-1}(q_3)]$, which we rewrite as $z_3 = y_1 \cap
f^{-1}[v \cap g^{-1}(q_3)]$.  But $v \cap g^{-1}(q_3)$ is a
finite collection of points.  If we add them with multiplicity we
get $L(g)$ points.  On the other hand, for each such point $p$ in
$S^{n-k_1+1}$, $y_1 \cap f^{-1}(p)$ is a finite collection of
points with total multiplicity $L(f)$.  Therefore, $z_3$ is a
finite collection of points with total multiplicity $L_2(F) =
L(f) L(g)$.

\vskip5pt

In the rest of this section, we prove the two propositions from
the introduction. These propositions are weaker than the main
theorem, and we include their proofs mostly for background.  The
main point of the paper is the improvement between Proposition 1.3
and Theorem 1.  The proofs of these propositions are easy
variations on the materical in \cite{G1}.  The proof of Proposition
1.2 is essentially due to Gromov.

\begin{prop} Suppose that $2 \le k_1 \le k_2 \le k_3$ and $n =
k_1 + k_2 + k_3 - 2$.  Let $F$ be a map from $(S^n, g)$ to
$(S^{k_1} \vee S^{k_2} \vee S^{k_3}, h_1 \vee h_2 \vee h_3)$ with
$k_1$-dilation at most $L^{k_1}$.  (For instance, $F$ may have
Lipschitz constant $L$.)  Then $L_2(F)$ is bounded as follows.

$$|L_2(F)| \le Iso_{n-k_1}(g) Iso_{k_3 - 1}(g) Vol(g) L^{n+2}
Vol(h_1)^{-1} Vol(h_2)^{-1} Vol(h_3)^{-1}. $$

\end{prop}

\proof By the coarea formula, we can choose 
$q_1$ so that $z_1 = F^{-1}(q_1)$ has volume at
most $L^{k_1} Vol(g) Vol(h_1)^{-1}$.  (See \cite{G1} for
more details.)

Then we can choose $y_1$ with volume at most
$Iso_{n-k_1}(g) L^{k_1} Vol(g) Vol(h_1)^{-1}$.

Using the coarea formula again, we choose $q_2 \in S^{k_2}$ so
that the volume of $z_2 = y_1 \cap F^{-1}(q_2)$ is at most
$Iso_{n-k_1} Vol(g) L^{k_1 + k_2} Vol(h_1)^{-1} Vol(h_2)^{-1}$.

Then we can choose $y_2$ with volume at most $Iso_{n-k_1}(g)
Iso_{k_3 -1 }(g) Vol(g) L^{k_1 + k_2}$ $Vol(h_1)^{-1}
Vol(h_2)^{-1}$.

But the degree of $\pi \circ F$ on $y_2$ is at most
$L^{k_3} Vol(y_2) Vol(h_3)^{-1}$.  Filling in our
bound for the volume of $y_2$ finishes the proof. \endproof

In \cite{G1}, we estimated the isoperimetric constants
of ellipses.

\begin{refprop} (\cite{G1}) Let $E$ be an n-dimensional
ellipse with principal axes $E_0 \le ... \le E_n$.  
Up to a constant factor $C(n)$,
$Iso_k(E) \sim E_{k+1} + E_{n-k}$.
\end{refprop}

Plugging this estimate into the last proposition, we immediately
get the following estimate.

\begin{prop} Suppose that $2 \le k_1 \le k_2 \le k_3$ and that $n +
2 = k_1 + k_2 + k_3$.  Let $E$ be an n-dimensional ellipse
with principal axes $E_0 \le ...
\le E_n$.  Let $F$ be a map from $E$ to the wedge of unit spheres
$S^{k_1} \vee S^{k_2} \vee S^{k_3}$ with $k_1$-dilation at most
$L^{k_1}$.  

If $k_3 < (n+1)/2$, then $L_2(F)$ is bounded as follows.

$$|L_2(F)| \le C(n) E_{n-k_1 + 1} E_{n-k_3 + 1} Vol(E) L^{n+2}.
$$

If $k_3 \ge (n+1)/2$, then $L_2(F)$ is bounded as follows.

$$|L_2(F)| \le C(n) E_{n-k_1 + 1} E_{k_3} Vol(E) L^{n+2}.
$$

\end{prop}

In this paper, we will study how sharp this estimate is and
improve it in the first case, $k_3 < (n+1)/2$.  In Section 5, we
will construct maps with large $L_2$ invariant, proving the
following theorem.

\begin{theorem 2} Suppose that $2 \le k_1 \le k_2 \le k_3$, and
$n + 2 = k_1 + k_2 + k_3$.  Let $E$ be an n-dimensional ellipse
with principal axes $E_0 \le ...
\le E_n$.  Suppose that $L > C(n) E_1^{-1}$.  Then there is a map
$F$ from $E$ to the wedge of unit spheres
$S^{k_1} \vee S^{k_2} \vee S^{k_3}$ with Lipschitz constant $L$
and $L_2(F)$ bounded below as follows.

$$L_2(F) \ge c(n) E_{n-k_1 + 1} E_{k_3} Vol(E) L^{n+2}.
$$
\end{theorem 2}

This theorem shows that our proposition is sharp up to a constant
factor in the case $k_3 \ge (n+1)/2$.  In the other case, it
turns out that the proposition is not sharp up to a constant
factor.  We will improve it in the next two sections.

\section{Directionally-dependent isoperimetric inequalities}

We begin by defining directional volume.  Let $C$ be an integral 
Lipschitz m-chain in $\mathbb{R}^n$.  Suppose that $J$ is an
m-tuple of distinct integers between 1 and $n$.  Let $P(J)$
denote the m-plane with coordinates $x_i$, $i \in J$. We define
the $J$-volume of $C$ to be the volume of the projection of $C$
to $P(J)$, counted with geometric multiplicity.  For example, if $J$ is any
(n-1)-tuple of numbers from 1 to n, and $C$ is the unit
(n-1)-sphere in $\mathbb{R}^n$, then the J-volume of $C$ is equal
to twice the volume of the unit (n-1)-ball.  

Here's another way of defining J-volume.  Let $TC_x$ denote the
tangent plane to $C$ at $x$.  For Lipschitz chains, $TC_x$ is
defined for almost every $x$ in $C$.  By an abuse of notation,
we write $TC_x \cdot P(J)$ to denote the inner product of
the unit k-vector corresponding to $TC_x$ and the unit k-vector
corresponding to $P(J)$.

$$Vol_J(C) := \int_C |TC_x \cdot P(J)| dvol(x).$$

The total volume of
$C$ is roughly equal to the sum of the volumes in different
directions.

$$Vol(C) \le \sum_J Vol_J(C) \le {n \choose m} Vol(C).$$

In \cite{LW}, Loomis and Whitney proved a directional estimate
for the volumes of open sets in $\mathbb{R}^n$.  Their original
estimate was written in terms of the projections of a set
to coordinate planes, but an immediate corollary is the following
estimate.

\begin{reftheorem} (Loomis, Whitney) Suppose that $H$ is a closed
embedded hypersurface in $\mathbb{R}^n$.  Let $V$ denote the volume
of the region enclosed by $H$.  This volume is bounded in terms
of the directional volumes of $H$ by the following formula.

$$V \le [\prod_J Vol_J(H)]^{\frac{1}{n-1}}.$$

Here the product is taken over the (n-1)-tuples $J$ of numbers from
1 to n.
\end{reftheorem}

We are interested in estimates that hold for cycles of any codimension.
The fundamental isoperimetric inequality for cycles of any codimension
was proven by Federer and Fleming.

\begin{reftheorem} (Federer, Fleming) Suppose that $z$ is a closed k-cycle
in $\mathbb{R}^n$.  Then there is a (k+1)-chain $y$ with $\partial y = z$
obeying the following volume bound.

$$|y| \le C(n) |z|^{\frac{k+1}{k}}.$$
\end{reftheorem}

There is a natural conjecture that generalizes the Loomis-Whitney theorem
to cycles of any codimension, which we include here for reference.

\begin{conj} Suppose that $z$ is a closed k-cycle in $\mathbb{R}^n$.
Then there is a (k+1)-chain $y$ with $\partial y = z$ so that for
every (k+1)-tuple $I$, the I-volume of $y$ is bounded in terms of the
directional volumes of $z$ as follows.

$$ Vol_I(y) \le C(n) [\prod_{J \subset I} Vol_J(z)]^{\frac{1}{k}}. $$

\end{conj}

(The product is taken over all k-tuples $J$ contained in $I$.  For each $I$,
there are (k+1) such k-tuples $J$.)

In this paper, we need directional isoperimetric estimates for cycles 
in an ellipse.  The estimate that we prove will depend
on the principal axes of the ellipse.  In fact, our goal is
to understand how the directional isoperimetric estimates depend on
the principal axes.

In the introduction to this paper, we mentioned an isoperimetric estimate 
for cycles in an ellipse from \cite{G1}.

\begin{iso} (\cite{G1}) Suppose that $E$ is an n-dimensional
ellipse with principal axes $E_0 \le ... \le E_n$.  Suppose that
$z$ is an integral m-cycle in $E$.  Then there is an (m+1)-chain
$y$ with $\partial y = z$ obeying the following estimate.

$$ |y| \le C(n) [E_{n-m} + E_{m+1}] |z|. $$

\end{iso}

We have defined the directional volumes for chains in Euclidean
space.  We extend the definition to chains in $E$ in the following
way.  The ellipse $E$ is $C(n)$-bilipschitz to the double of a
rectangle $R$ with dimensions $E_1 \le ... \le E_n$.  We fix
a particular bilipschitz equivalence.  Now given an m-chain $z$ in
the double of $R$, we let $z_N$ be the intersection of $z$ with
the Northern hemisphere and we let $z_S$ be the intersection of
$z$ with the Southern hemisphere.  We view $z_N$ and $z_S$ as
chains in the rectangle $R$, and so we know how to define their
directional volumes.  Then we define the $J$-volume of $z$ to be
$Vol_J(z_N) + Vol_J(z_S)$.

Now we refine the isoperimetric inequality above, taking into
account the directional volumes of $y$ and $z$.

\begin{diriso} Let $z$ be an integral m-cycle in $E$. 
Then $z$ bounds an (m+1)-chain $y$ with the following bounds on
directional volumes.  For each (m+1)-tuple $I$, we let $i$ denote
the smallest number in $I$ and $e$ denote the smallest number not
in $I$.

$$Vol_I(y) \le C(n) \left[ E_i Vol_{I-i}(z) + \sum_{d=1}^{e-1} E_d
Vol_{I-d} (z) \right]. $$

\end{diriso}

Since $i \le n-m$ and $d \le e-1 \le m+1$, we see that the total
volume of $y$ is bounded by $C(n) [E_{n-m} + E_{m+1}] Vol(z)$ recovering
the standard isoperimetric inequality in $E$.

Our directional isoperimetric inequality improves on the standard one
in two ways.  First, if we input a cycle $z$ with only a bound on the
total volume of $z$, then we get out a chain $y$ whose total volume
obeys the standard bound, but which has smaller directional volumes in
most directions.  Second, if we input a cycle $z$ with some control
on the directional volumes, then we may be able to output a chain $y$
with smaller total volume then the standard isoperimetric inequality
can deliver.

The proof of this directional isoperimetric inequality is a more 
complicated version of the proof of
the standard isoperimetric inequality in \cite{G1}.  
Since the proof below is somewhat involved,
it might help the reader to look at the proof in \cite{G1} first.  

We
build up to the result we need in three steps.  First we prove an
estimate for absolute cycles in a rectangle.  Second we prove an
estimate for relative cycles in a rectangle.  Third, we combine
these results to get an estimate for cycles in an ellipse.

\vskip5pt

{\bf Isoperimetric inequality for absolute cycles in a rectangle}

\begin{prop} Suppose that $z$ is an m-dimensional cycle in the
rectangle $R$ with dimensions $R_1 \le ... \le R_n$.
Then there is an (m+1)-chain $y$ in $R$ with $\partial y = z$
obeying the following bounds.

For each (m+1)-tuple $I$, let $i$ denote the smallest element
of $I$.  Let $I-i$ denote the m-tuple formed by removing $i$
from $I$.

$$Vol_I(y) \le R_i Vol_{I-i}(z).$$

\end{prop}

\proof We proceed by induction on the dimension $n$.  The result
is vacuous when the dimension is zero.

Let $\pi$ denote the projection from $\mathbb{R}^n$ onto the
plane $x_1 = 0$.  There is a homology from $z$ to $\pi(z)$
consisting of a union of lines.  (It lies in the cylinder $\pi(z)
\times [0,1]$.) This homology has $(1 \cup J)$-volume at most $R_1
Vol_J(z)$ for each $J$ that does not contain $1$.  It has
$I$-volume $0$ for any $I$ that does not contain $1$.

By induction, $\pi(z)$ bounds a chain $y' \subset \{ 0 \} \times
[0, R_2] \times ... \times [0, R_n]$, where the $I$-volume of
$y'$ is bounded by $R_i Vol_{I-i}(\pi(z)) \le R_i
Vol_{I-i}(z)$, where $I$ is any (m+1)-tuple of $2.. n$. 
Assembling the first homology with the filling $y'$ finishes the
proof. \endproof

Remark: For each $I$, we have bounded the I-volume of $y$ in
terms of only the $I-i$ volume of $z$.  If $J$ is an m-tuple
containing 1, then $J$ is not equal to $I-i$ for any (m+1)-tuple
$I$.  Therefore, we can bound the total volume of $y$ using only
some of the directional volumes of $z$.  We will need this observation
in the proof of our isoperimetric inequality for ellipses.

\vskip5pt

{\bf Isoperimetric inequality for relative cycles in a rectangle}

\vskip5pt

Next we study relative cycles in a rectangle.  We think of a
relative cycle as a chain $z$ with $\partial z$ contained in the
boundary of the rectangle.  We would like to ``push'' $z$ into
the boundary.  In other words, we want to find an (m+1)-chain $y$
with $\partial y = z + B$ where $B$ is contained in $\partial R$. 
For our purposes, we need
estimates for both the size of $y$ and the size of $B$.

\begin{prop} Suppose that $z$ is an m-dimensional relative
integral cycle in the rectangle $R$.  Then there is an
(m+1)-chain $y$ with $\partial y = z + B$ and $B$ contained in
$\partial R$ obeying the following estimates.

Let $I$ be an (m+1)-tuple.  Let $e$ denote the smallest number not in $I$.

$$Vol_I(y) \lesssim \sum_{d=1}^{e-1} R_d Vol_{I-d}(z). \eqno{(*)}$$

If $J$ is an m-tuple that does not include $1$, and if
$e$ denotes the smallest number in $[2.. n]$ which does not lie in $J$,
then the $J$-volume of $B$ obeys the following estimate.

$$Vol_J(B) \lesssim Vol_J(z) + \sum_{d=1}^{e-1} R_1^{-1} R_d
Vol_{1 \cup J - d}(z). \eqno{(**)}$$  

\end{prop}

(We are only able to bound some of the directional volumes of
$B$.  If $J$ includes $1$, then we do not prove any upper
bound on $Vol_J(B)$.)

\proof We begin by proving a lemma that covers a special case. 
The special case occurs when the boundary of $z$ lies only in the
bottom and sides of $\partial R$ and does not touch the top of
$\partial R$.  The following figure illustrates an example of a
relative cycle $z$ in this special case.

\includegraphics{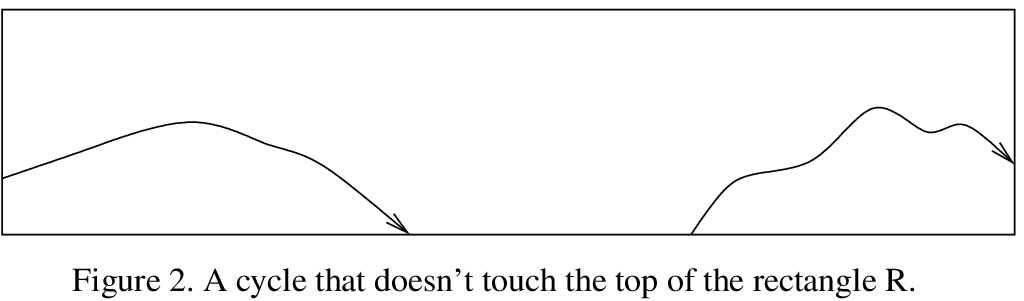}

\begin{lemma} Suppose that $z$ is an m-chain in $R$ with
$\partial z$ lying in $\partial R$.  Let $R'$ be the (n-1)-dimensional
rectangle $[0, R_2] \times ... \times [0, R_n]$ so that $R = [0, R_1]
\times R'$.  Suppose that the boundary of $z$ does not intersect
$\{ R_1 \} \times R'$.  Then there is an (m+1)-chain $y$ with
$\partial y = z + B$ where $B$ is an m-chain in $\partial R$
obeying the following inequalities.

1. If $I$ is an (m+1)-tuple containing $1$, then $Vol_I(y) \le
R_1 Vol_{I-1}(z)$.

2. If $I$ is an (m+1)-tuple that does not contain $1$, then
$Vol_I(y) = 0$.

3. If $J$ is an m-tuple that does not contain $1$, then $Vol_J(B)
\le Vol_J(z)$.

\end{lemma}

\proof Let $I: z \rightarrow R$ be
the identity embedding, and let $I_1, ... I_n$ be its $n$
coordinates.  Now we construct a map $f: z \times
[0,1] \rightarrow R$ with coordinate functions defined as
follows: $f_1(x,t) = t I_1(x)$ and for all $i \not= 1$, $f_i(x,t) =
I_i(x)$. We define the filling $y$ to be $f(z \times [0,1])$.

\includegraphics{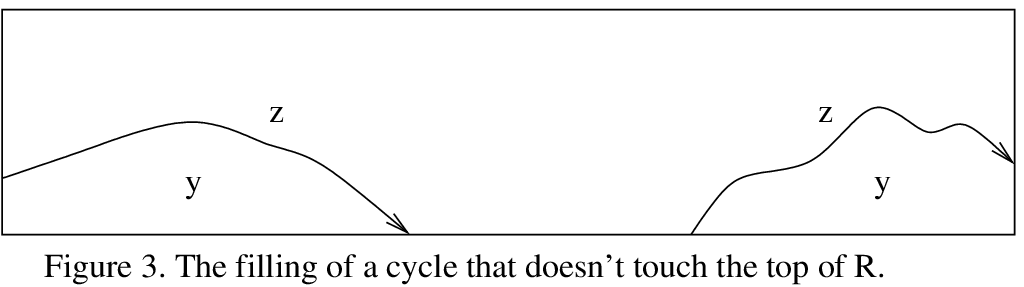}

First we check that $y$ is a filling of $z$.  The boundary of $y$
is equal to $f(\partial z \times [0,1]) - f(z \times \{0\})
+ f(z \times \{ 1 \})$. Let $x$ be a point in $\partial z$. 
If $x$ lies in a "side" of $R$ (i.e. in $[0,R_1] \times
\partial R'$), then $f(x,t)$ lies in that "side" of $R$ for all
$t$.  If $x$ lies in the "bottom" of $R$ (i.e. in $\{ 0 \} \times
R'$), then $f(x,t)$ lies in the bottom of $R$ for all $t$.  By
assumption, $\partial z$ lies in the sides and bottom of $R$.
Therefore, $f(\partial z \times [0,1])$ lies in the boundary of
$R$.  Also $f(z \times \{ 0 \})$ lies in the bottom of $R$. 
Finally $f$ restricted to $z \times \{ 1 \}$ is the identity, and
so $f(z \times \{ 1 \})$ is just $z$.  Hence the boundary of $y$
is equal to $z$ plus a chain lying in $\partial R$.  We call this
chain $B$.

Now it remains to prove our estimates for $y$ and $B$.  Let $\pi$
denote the projection from $R$ onto $R'$.  Then $y$ lies in $[0,
R_1] \times \pi(z)$. 
That proves estimates 1 and 2.  According to the last paragraph,
$B$ is made up of two pieces: $- f(\partial z \times [0,1])$ and
$- f(z \times \{ 0 \})$.  The first piece has $J$-volume equal to
zero unless $1$ lies in $J$.  The last piece is just the
projection $\pi(z)$.  Since $Vol_J(\pi(z)) \le Vol_J(z)$ for any
$J$, we get the last inequality. \endproof

To prove Proposition 2.2, we reduce our situation to the case of
the lemma by cutting an arbitrary cycle $z$ into pieces in such
a way that each piece can be filled either by using Lemma 2.1 or
by induction on the dimension.  By induction, we can assume that
Proposition 2.2 holds for rectangles of dimension n-1.

We consider the slices $z_h := z \cap \{ x | x_1 = h \}$ for
various heights $h$.  By the coarea inequality, we can choose $h$
so that the following inequality holds for each (m-1)-tuple $K$
in $2... n$.

$$Vol_K (z_h) \lesssim R_1^{-1} Vol_{1 \cup K} (z).$$

(On the other hand, if $K$ contains 1, then the K-volume of
$z_h$ is zero.)

The following picture shows an example of $z_h$.

\vskip5pt

\includegraphics{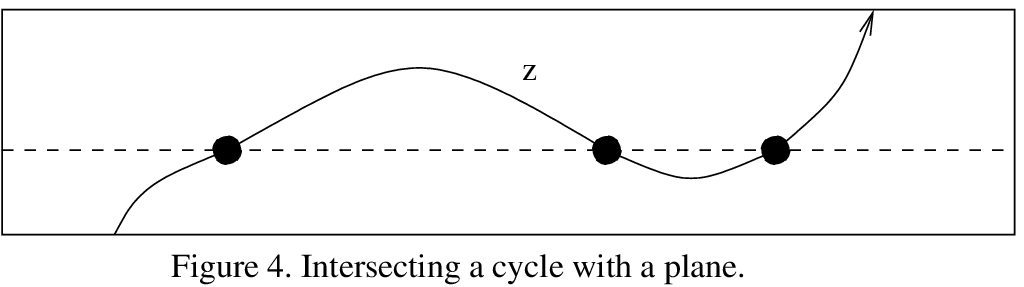}

In this figure, the solid oriented curve denotes the cycle $z$. 
The dotted line denotes the plane $x_1 = h$.  The three dark
points denote their intersection, $z_h$.

We now decompose $z$ into two pieces as follows.

$$z = \bigl( z - [0,R_1] \times z_h \bigr) + [0,R_1] \times z_h =
z_1 + z_2.$$

We deal with the first piece by decomposing it into an upper and lower
half: $z_1 = z_+ + z_-$ where $z_+$ is the part of $z_1$ lying
above $x_1=h$ and $z_-$ is the part of $z_1$ lying below $x_1=h$. 
The chains $z_+$ and $z_-$ are each relative cycles in $R$.  The
following figure shows $z_+$ and $z_-$ in the example from the
last figure.

\vskip5pt

\includegraphics{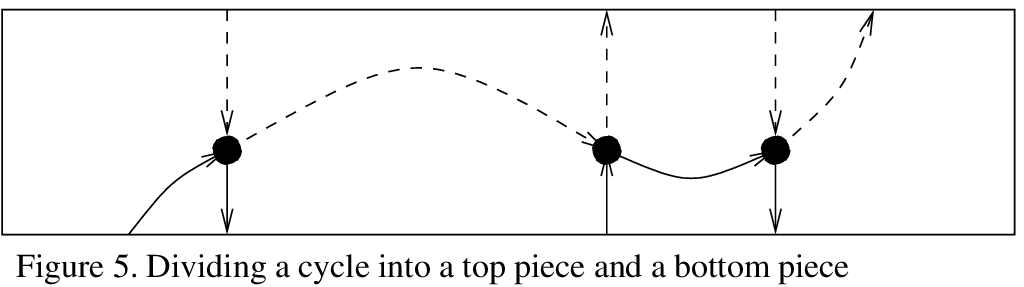}

In this figure, the dotted curves denote the cycle $z_+$ and the
solid curves denote the cycle $z_-$.  As in the previous figure,
the three dark points denote $z_h$.

The cycle $z^+$ avoids the bottom of $R$ and $z^-$ avoids the top
of $R$, and so we can fill them both using Lemma 2.1.  Let $y_+$
be the filling of $z_+$ and $y_-$ the filling of $z_-$. 
According to Lemma 2.1, the directional volumes of $y_+$ are
bounded in terms of the directional volumes of $z^+$.  But
$Vol_J(z^+) \le Vol_J(z) + Vol_J(z_h \times [0,R_1])$.  We chose
$h$ so that $Vol_K(z_h) \lesssim R_1^{-1} Vol_{1 \cup K} (z)$ for
each (m-1)-tuple $K$, and so $Vol_J(z_h \times [0, R_1]) \lesssim
Vol_J(z)$ for every $J$.  Then the conclusion of Lemma 2.1 shows
that $y_+$ obeys $(*)$ and $B_+$ obeys $(**)$.  The same holds
for $z_-$ and its filling $y_-$.  Finally we define $y_1 = y_+ +
y_-$ and $B_1 = B_+ + B_-$.  We have seen that $\partial y_1 =
z_1 + B_1$, and that $B_1$ lies in the boundary of $R$, and that
$y_1$ and $B_1$ obey the directional volume estimates $(*)$ and
$(**)$.

Now we have reduced matters to a cycle $z_2$ of the special form
$[0, R_1] \times z_h$.  By induction on the dimension, 
there is an m-chain $y_h$ in $[0,R_2] \times ...
\times [0, R_n]$ with $\partial y_h = z_h + B_h$ obeying the
estimate $(*)$.  If we spell out $(*)$ we get the following.

Let $J$ be any m-tuple in $[2.. n]$.  Suppose that $e$ denotes
the smallest number in $[2.. n]$ but not in $J$.

$$Vol_J(y_h) \lesssim \sum_{d=2}^{e-1} R_d Vol_{J-d}(z_h) \lesssim 
\sum_{d=2}^{e-1} R_1^{-1} R_d Vol_{1 \cup J - d}(z). \eqno{(1)}$$

Now we define $y_2 = [0, R_1] \times y_h$.  
Suppose that $I$ includes 1.  Let $J$ denote
$I-1$.  Then let $e$ denote the smallest number not in $I$, which
is the same as the smallest number in $[2.. n]$ but not in $J$.

$$Vol_I(y_2) = R_1 Vol_J(y_h) \lesssim \sum_{d=2}^{e-1} R_d
Vol_{I-d}(z).$$

On the other hand, if $I$ does not include 1, then the $I$-volume
of $y_2$ is zero.  So $y_2$ obeys inequality $(*)$.

Next we define $B_2$ by setting $\partial y_2 = z_2 + B_2$. 
Since $y_2 = [0,R_1] \times y_h$, the boundary $\partial y_2$ is
equal to $[0, R_1] \times z_h + [0, R_1] \times B_h - \{ 0 \}
\times y_h + \{ 1 \} \times y_h$.  The first term $[0, R_1]
\times z_h$ is $z_2$, and so the remaining terms are equal to
$B_2$.  The three terms making up $B_2$ each lie in $\partial R$.

Finally, we have to bound (some of) the directional volumes of
$B_2$.  Suppose that $J$ is an m-tuple which does not contain
$1$.  The J-volume of $[0, R_1] \times B_h$ is zero.  Each of the
other two terms has J-volume equal to that of $y_h$.  Applying
$(1)$, we get the following estimate for any m-tuple $J$ which
does not contain 1.

$$Vol_J(B_2) \lesssim \sum_{d=2}^{e-1} R_1^{-1} R_d Vol_{1 \cup J
- d}(z).$$

This equation shows that $B_2$ obeys $(**)$.

Finally, we set $y = y_1 + y_2$ and $B = B_1 + B_2$.  Now we have
$\partial y = z + B$, where $B$ lies in $\partial R$, and $y$ and
$B$ obey inequalities $(*)$ and $(**)$. \endproof

\vskip5pt

{\bf Isoperimetric inequality for cycles in an ellipse}

\vskip5pt

Let $E$ be the n-dimensional ellipse with principal axes 
$E_0 \le ... \le E_n$.  We recall the definition of the
directional volume for chains in $E$.  The ellipse $E$
is $C(n)$-bilipschitz to the double of a rectangle $R = [0, E_1]
\times ... \times [0, E_n]$.  We fix a bilipschitz equivalence to
use throughout the paper.  Each copy of $R$ in the bilipschitz
equivalence is a hemisphere of $E$.  Now given an m-chain $z$ in
the double of $R$, we let $z_N$ be the intersection of $z$ with
the Northern hemisphere and we let $z_S$ be the intersection of
$z$ with the Southern hemisphere.  We view $z_N$ and $z_S$ as
chains in the rectangle $R$, and so we know how to define their
directional volumes.  If $J$ is any m-tuple of the numbers from
1 to n, we define the $J$-volume of $z$ to be $Vol_J(z_N) + Vol_J(z_S)$.

With this definition of directional volume, we can state our
directional isoperimetric inequality.

\begin{prop} Let $z$ be an m-dimensional integral cycle in $E$. 
Then $z$ bounds an (m+1)-chain $y$ with the following bounds on
directional volumes.  For each (m+1)-tuple $I$, we let $i$ denote
the smallest number in $I$ and $e$ denote the smallest number not
in $I$.

$$Vol_I(y) \lesssim E_i Vol_{I-i}(z) + \sum_{d=1}^{e-1} E_d
Vol_{I-d} (z). $$

\end{prop}

\proof This proposition follows by combining the previous two.
As above, we let $z_S$ be the intersection of $z$ with the
Southern hemisphere.  The chain $z_S$ is a relative cycle.  We
apply Proposition 2.2, which tells us that there is a chain $y_S$ in
the Southern hemisphere with $\partial y_S = z_S + B$ and $B
\subset \partial R$ obeying the following estimates.

1. Let $I$ be an (m+1)-tuple and $e$ the smallest number not in
$I$. Then $Vol_I(y_S) \lesssim \sum_{d=1}^{e-1} E_d
Vol_{I-d}(z_S)$.

2. Let $J$ be an m-tuple not containing $1$ and $f$ the smallest
number in $[2.. n]$ not in $J$.  Then $Vol_J(B) \lesssim
Vol_J(z_S) + \sum_{d=1}^{f-1} E_1^{-1} E_d Vol_{1 \cup J - d}(z_S)$.

Now $z_N - B$ is an absolute m-cycle in the Northern hemisphere. 
We apply the directional isoperimetric inequality for absolute
cycles to fill it.  This inequality tells us that there is a
chain $y_N$ in the Northern hemisphere with $\partial y_N = z_N -
B$ obeying the following estimate.

3. Let $I$ be an (m+1)-tuple and let $i$ be the smallest number in $I$.
Then $Vol_I(y_N) \lesssim E_i [Vol_{I-i}(z_N) + Vol_{I-i}(B)]$.

Since $i$ is the smallest number in $I$, $I-i$ does not contain $1$, and
so we can use 2 to bound $Vol_{I-i}(B)$.  We do this in two
cases.  First we consider the case $i > 1$.  In this case $I-i$
does not contain 2.  Hence $f=2$ in inequality 2, and we conclude
that $Vol_{I-i}(B) \lesssim Vol_{I-i}(z_S)$.  Plugging this inequality
into 3, we get the following estimate.

4a. If $i > 1$, then $Vol_I(y_N) \lesssim E_i Vol_{I-i}(z)$.

Next we consider the case $i=1$.  We recall that $e$ is the
smallest number not in $I$.  Hence $e$ is also the smallest
number in $[2.. n]$ which is not in $I-1$.  In this case,
inequality 2 tells us that $Vol_{I-1}(B) \lesssim E_1^{-1}
\sum_{d=1}^{e-1} E_d Vol_{I-d}(z)$.  Plugging this estimate into
inequality 3, we get the following.

4b. If $i=1$, then $Vol_I(y_N) \lesssim \sum_{d=1}^{e-1} E_d
Vol_{I-d}(z)$.

Finally, we let $y = y_N + y_S$.  The boundary $\partial y = z$. 
Combining estimates 1, 4a, and 4b, we see that the directional
volumes of $y$ obey the conclusion of the proposition. \endproof

\section{A coarea inequality for directional volumes}

In order to bound the second-order linking invariant of a
Lipschitz map, we also need a directional version of the coarea
inequality.  This inequality is only a minor variation on the
standard one.  The proof combines the general coarea formula with
some calculations in exterior algebra.

\begin{prop} Let $y$ be a $C^\infty$ m-chain in $\mathbb{R}^n$, and let
$F$ be a $C^\infty$ map from $y$ to $(N^q, h)$ with q-dilation at
most $\Lambda$.  Then $F$ has a fiber $z = F^{-1}(n)$ for some $n
\in N$ obeying the following estimates for the directional
volumes.  Let $k= m-q$ be the dimension of $z$ and let $J$ be any
k-tuple

$$Vol_J(z) \le C(n) \Lambda Vol(N)^{-1} \sum_{J \subset I}
Vol_I(y).$$

\end{prop}

\proof First we write down the general coarea formula, which
holds for any function $G$ on $y$.

$$\int_y Jac[dF(x)] G(x) dvol(x) = \int_N [\int_{F^{-1}(n)} G
dvol_{F^{-1}(n)}] dvol_h(n).$$

At points where $Jac[dF(x)] \not= 0$, the kernel of $dF$ is a
k-plane.  We write $V(F)$ to denote the unit k-vector parallel to
this k-plane.  (We should specify a choice of orientation, but
the orientations won't matter because we will always take
absolute values.)  Then we take $G = |J \cdot V(F)|$.  With this
choice, the integral over the fiber $\int_{F^{-1}(n)} G$ is
exactly the J-volume of $F^{-1}(n)$.  Therefore, we have the
following formula.

$$\int_N Vol_J[F^{-1}(n)] dvol_h(n) = \int_y Jac[dF(x)] |J \cdot
V(F)| dvol \le \Lambda \int_y |J \cdot V(F)| dvol.$$

We don't know in which direction the plane $V(F)$ points, except
that it is a subplane of the tangent plane to $y$.  Let $P$
denote the tangent space to $y$ at a given point.  We are led to
estimate $\sup_{Q^k \subset P^m} |J \cdot Q|$, the largest
possible value of the term $|J \cdot V(F)|$.  This is a problem about
exterior algebra which turns out to have a clean answer.

\begin{lemma} Let $P^m$ denote an m-dimensional plane in $\mathbb{R}^n$.
Let $J$ be a k-tuple.  By abuse of notation, we also let 
$J$ denote the unit k-vector corresponding
to the k-tuple $J$.

$$\sup_{Q^k \subset P^m} |J \cdot Q| = \left[ \sum_{J \subset I}
|I \cdot P|^2 \right]^{1/2}.$$

\end{lemma}

\proof Both sides are invariant if we rotate the plane $P$ that leaves the plane
spanned by $J$ invariant.  By using
such a rotation, we can arrange that $P \cap J^\perp$ is in
standard position.  We let $K$ be an (m-k)-tuple disjoint from
$J$.  Because of the rotational symmetry, we can assume without
loss of generality that $P \cap J^\perp = K$.  We let $Q_0
= P \cap K^{\perp}$.  Hence $P = Q_0 \oplus K$.

On the one hand, $\sup_{Q^k \subset P^m} |J \cdot Q| = |J \cdot
Q_0|$.  To see this, let $Q$ be any plane in $P$ and write its
fundamental k-vector as a wedge of unit vectors $v_1 \wedge ...
\wedge v_k$ with $v_i$ in $P$.  Decompose each vector $v_i$ into
a piece in $Q_0$ and a piece in $K$, $v_i = u_i + w_i$. 
Expanding the wedge product, we get a sum of terms.  Each term
involving any $w_i$ vanishes when we take the inner product with
$K$.  The other term is equal to $c Q_0$ for some $c$ with $|c|
\le 1$.

On the other hand, the right-hand side is also equal to $|J \cdot
Q_0|$.  The right-hand side vanishes unless $I = J \cup K$, and
so the right hand side is $|J \wedge K \cdot Q_0 \wedge K| = |J
\cdot Q_0|$.  \endproof

Applying Lemma 3.1 we get the following.

$$\int_N Vol_J[F^{-1}(n)] dvol_h(n) \le \Lambda \int_{y} [\sum_{J
\subset I} |I \cdot Ty|^2]^{1/2} dvol
\le \Lambda \sum_{J \subset I} Vol_I(y).$$

This formula holds for each choice of $J$.  Therefore, we may
choose $n \in N$ so that for every k-tuple $J$, the following
holds.

$$ Vol_J(F^{-1}(n)) \le {n \choose k} \Lambda Vol(N)^{-1} 
\sum_{J \subset I} Vol_I(y).$$

\endproof

\section{Dilations and second-order linking invariants}

In this section we prove Theorem 1.

\begin{theorem 1} Suppose that $2 \le k_1 \le k_2 \le k_3$, that $n +
2 = k_1 + k_2 + k_3$.  Let $E$ be an n-dimensional ellipse 
with principal axes $E_0 \le ... \le E_n$.  Let $F$ be a
map from $E$ to the wedge of unit spheres $S^{k_1} \vee S^{k_2}
\vee S^{k_3}$ with $k_1$-dilation at most $L^{k_1}$.  Then
$L_2(F)$ is bounded as follows.

$$|L_2(F)| \le C(n) E_{n-k_1 + 1} E_{k_3} Vol(E) L^{n+2}.
$$

\end{theorem 1}

The idea of the proof is to imitate the argument in the proof
of Proposition 1.3 but to substitute the directional isoperimetric
for the standard isoperimetric inequality.

\proof We begin by choosing a point $q_1$ in $S^{k_1}$ and looking
at the fiber $z_1 = F^{-1}(q_1)$.  For generic $q_1$, the inverse
image is a manifold of dimension $n-k_1$.  By the coarea formula,
we can choose $q_1$ so that the fiber $F^{-1}(q_1)$ has volume at
most $Vol(E) L^{k_1} / Vol(S^{k_1}) \le C(n) Vol(E) L^{k_1}$.  In
particular, each $J$-volume of $z_1 = F^{-1}(q_1)$ is at most
$C(n) Vol(E) L^{k_1}$.

Next we choose a chain $y_1$ with boundary $z_1$, using the
directional isoperimetric inequality Proposition 2.3.  The chain
$y_1$ will obey the following directional volume bounds.

Let $I$ be a $(n-k_1+1)$-tuple.  Let $i$ denote the smallest
element in $I$ and let $e$ denote the smallest element not in $I$.

$$Vol_I(y_1) \le C(n) Vol(E) L^{k_1} [E_i + E_{e-1}]. \eqno{(1)}$$

Next we choose a point $q_2$ in $S^{k_2}$ and look at the
intersection $z_2 = y_1 \cap F^{-1}(q_2)$. 
By using the directional coarea inequality, we can bound the
directional volumes of $z_2$.

The cycle $z_2$ has dimension $n-k_1+1-k_2$.  Let $J$ be a tuple
of that dimension, and let $I$ be an $(n-k_1+1)$-tuple containing
$J$.  Let $i$ be the smallest element in $I$ and let $e$ be the
smallest element not in $I$.  Because of the cardinality of $I$,
$i \le k_1$.  The tuple $I$ is formed by adding $k_2$ elements to
the tuple $J$.  Let $f$ denote the $(k_2+1)^{st}$ smallest element
which is not in $J$.  (In other words, we list the elements not in
$J$ from smallest to largest, and let $f$ be the $(k_2 + 1)^{st}$
element in this list.)  Then $e
\le f$.  Therefore we get the following estimate for the
$J$-volumes of $z_2$.

$$Vol_J(z_2) \le C(n) Vol(E) L^{k_1+k_2} [E_{k_1} + E_{f-1}].
\eqno{(2)}$$

The third step of the proof is to apply the directional
isoperimetric inequality again to estimate the size of a filling
$y_2$ of $z_2$.

The chain $y_2$ has dimension $n- k_1 - k_2 + 2 = k_3$.  Let $K$
denote a tuple of that dimension.  Let $k$ denote the smallest
element in $K$ and let $g$ denote the smallest element not in
$K$.  Proposition 2.3 gives the following estimate.

$$Vol_K(y_2) \le C(n) [E_k Vol_{K-k}(z_2) + \sum_{c=1}^{g-1} E_c
Vol_{K-c}(z_2)].$$

Next we plug in the estimate for the $J$-volume of $z_2$ from
equation 2.  We use this estimate to substitute for $Vol_{K-k}(z_2)$
and $Vol_{K-c}(z_2)$.

$$Vol_K(y_2) \lesssim Vol(E) L^{k_1 + k_2} \left[ E_k (E_{k_1} + E_{f(K-k)-1})
+ \sum_{c=1}^{g-1} E_c (E_{k_1} + E_{f(K-c)-1} ) \right].$$

We will check that the bracketed expression is bounded by
$E_{n-k_1+1} E_{k_3}$.

First we deal with the term $E_k E_{k_1}$.  Recall that $k$ is
the smallest element in $K$.  Because the cardinality of $K$ is
$k_3$, $k \le n-k_3+1$.  Hence $E_{k} E_{k_1} \le E_{n-k_3+1}
E_{k_1} \le E_{n-k_1+1} E_{k_3}$. 

Second we deal with the term $E_k E_{f(K-k) - 1}$. Recall that
$f(K-k)$ is the $(k_2 + 1)^{st}$ smallest element not in $K-k$. 
The cardinality of $K-k$ is $n- k_1 -k_2 + 1$, and so $f \le n-k_1
+ 2$. Hence if $k \le k_3$, then $E_k E_{f(K-k) - 1} \le
E_{n-k_1+1} E_{k_3}$.  On the other hand, if $k \ge k_3+1 \ge k_2
+ 1$, then the numbers $1, ..., k_2+1$ are all not in $K-k$, and
so $f(K-k) = k_2+1$.  In this case $E_k E_{f(K-k) -1} \le
E_{n-k_3+1} E_{k_2} \le E_{n-k_1+1} E_{k_3}$.

Third we deal with the term $E_{k_1} E_c$.  
The cardinality of $K$ is $k_3$ and $g$ is
the smallest element not in $K$, so $g \le k_3 + 1$.  Since $c
\le g-1$, it follows that $c \le k_3$.  Hence $E_c E_{k_1}$ is
bounded by $E_{k_3} E_{k_1}$ which is bounded by $E_{n-k_1+1}
E_{k_3}$.

Finally, we deal with the term $E_c E_{f(K-c) - 1}$.  Recall 
that $f(K-c)$ is the $(k_2+1)^{st}$ smallest element not in
$K-c$.  Since $K-c$ has $k_3-1$ elements, $f(K-c) \le k_2 + k_3 =
n-k_1+2$.  Therefore, $E_{f(K-c) - 1} \le E_{n-k_1+1}$.  In the
last paragraph, we saw that $c \le k_3$.  So the product $E_c E_{f(K-c)-1}
\le E_{n-k_1+1} E_{k_3}$.  

Putting together the different
terms, we have bounded the total volume of $y_2$ as follows.

$$Vol(y_2) \le C(n) L^{k_1+k_2} Vol(E) E_{n-k_1+1} E_{k_3}.$$

Finally, we choose a point $q_3$ in $S^{k_3}$ and look at the
fiber $z_3 = y_2 \cap F^{-1}(q_3)$.  By the coarea inequality, we
can choose $q_3$ so that the 0-dimensional volume of $z_3$ is
at most $C(n) L^{n+2} Vol(E) E_{n-k_1+1} E_{k_3}$.  Since
$L_2(F)$ is the sum of the points in $z_3$ counted with
multiplicity $\pm 1$, we see that the norm of $L_2(F)$ is bounded
by $C(n) L^{n+2} Vol(E) E_{n-k_1+1} E_{k_3}$. \endproof

\section{Lipschitz maps with large linking invariants}

In this section we prove Theorem 2.

\begin{theorem 2} Suppose that $2 \le k_1 \le k_2 \le k_3$, that $n +
2 = k_1 + k_2 + k_3$.  Let $E$ be an n-dimensional ellipse with 
principal axes $E_0 \le ...
\le E_n$.  Suppose that $L > C(n) E_1^{-1}$.  Then there is a map
$\Phi$ from $E$ to the wedge of unit spheres
$S^{k_1} \vee S^{k_2} \vee S^{k_3}$ with Lipschitz constant $L$
and $L_2(\Phi)$ bounded below as follows.

$$L_2(\Phi) \ge c(n) E_{n-k_1 + 1} E_{k_3} Vol(E) L^{n+2}.
$$

\end{theorem 2}

\proof We begin by constructing a map with a large linking
invariant.  During the proof, we will use the map twice.  In
order for the proof to fit together, we need to carefully choose
the range of the map.

We will use the following vocabulary.  For any dimension $d \le
n$, we let $E[d]$ denote the d-dimensional ellipse with
principal axes $E_0 \le ... \le E_d$.  Note that $E[d]$ is
$C(n)$-bilipschitz to the double of the rectangle $[0, E_1]
\times ... \times [0, E_d]$.

The domain of the map is the ellipse $E$.  We will have to keep
track of the ``tips'' of the ellipse $E$.  If $E$ is given by
the equation $\sum_{i=0}^n (x_i/E_i)^2 = 1$, then the tips of $E$
are the two points $(0, ..., 0, \pm E_n)$.

The topology of the range is as follows.  Let $d,e$ be integers
at least 2 so that $n+1 = d + e$.  Let $p, q$ be antipodal points
on the sphere $S^d$.  Let $*$ be a basepoint of the sphere $S^e$. 
The range of our map is the space $X$ given by taking the union
of $S^d$ and $S^e$ and then identifying the points $p, q,$ and
$*$.  This identified point is the basepoint of $X$.

Next we define a metric on the space $X$, which just means
picking a metric on $S^d$ and a metric on $S^e$.  The metric on
$S^d$ is the ellipsoidal metric $E[d]$.  The two points $p,q$ are
the tips of the ellipse.  (If the ellipse is given by the
equation $\sum_{i=0}^d (x_i/E_i)^2 =1 $ in $\mathbb{R}^{d+1}$,
then the points $p,q$ are the points $(0, ..., 0, \pm E_d)$.)

The metric on $S^e$ is the one-point compactification of a
rectangle with dimensions $E_{n-e+1} \times ... \times E_n$.  In
other words, the metric is given by taking the Euclidean
rectangle $[0, E_{n-e+1}] \times ... \times [0, E_n]$ and
collapsing the boundary to a point.  The basepoint $*$ is the
point we added to do the compactification, or in other words the
point corresponding to the boundary.  (Our metric is singular at
the base point, but the singularity does not create any
problems.)  We write $(X,h)$ to refer to the space $X$ equipped
with this metric.

A map $F: S^n \rightarrow X$ has a linking invariant defined in
the same way as for a map from $S^n$ to $S^d \vee S^e$.  Namely,
let $q_1$ be a generic point of $S^d \subset X$ and $q_2$ a
generic point of $S^e \subset X$, and look at the linking number
of the two disjoint cycles $F^{-1}(q_1)$ and $F^{-1}(q_2)$.

\begin{lemma} There is a map $F: E \rightarrow (X,h)$ with
linking invariant 1 and with Lipschitz constant at most $C(n)$. 
Moreover, this map takes the tips of $E$ to the basepoint of $X$.
\end{lemma}

\proof We begin by writing down two open sets inside of $E$. 
Geometrically, the open sets are thick linked spheres.  In order
to write them down, we think of $E$ as the double of the
rectangle with dimensions $E_1 \times ... \times E_n$.

The first set $U$ has the form $S^{n-e} \times B^{e}$.  It is
the double of the following product: $\prod_{i=1}^{n-e} [0,E_i]
\times \prod_{i=n-e+1}^n [(1/3) E_i, (2/3) E_i]$.  A core sphere
$S^{n-e}$ in $U$ is given by the double of $\prod_{i=1}^{n-e} [0,
E_i]$ times a point. 

We define the map $F$ on $U$ as follows.  We let $\pi_U$
denote the projection from $U$ to the rectangle $R(U) =
\prod_{i=n-e+1}^n [(1/3) E_i, (2/3) E_i]$.  There is a degree one
map $\psi_U$ from $R(U)$ to $(S^e, h)$ taking the boundary of $R(U)$
to the basepoint and with Lipschitz constant 3.  The map $F$
on $U$ is given by the composition $\psi_U \circ \pi_U$.  It maps the
boundary of $U$ to the basepoint of $S^e$.  Since $S^e \subset
X$, we can think of $F$ as a map from $(U, \partial U)$ to $(X,
*)$.

We let $V_0$ be the
rectangle $\prod_{i=n-e+1}^n [(1/10) E_i , (9/10) E_i]$ minus the
interior of $R(U)$.  The set $V_0$ is homeomorphic to $S^{e-1}
\times [0,1]$.  Now the set $V$ is the double of the product
$\prod_{i=1}^{n-e} [0, E_i] \times V_0$.  Therefore $V$ is
homeomorphic to $S^{e - 1} \times S^{n-e} \times [0,1]$.  We call
a copy of $S^{e - 1}$ times a point a core sphere of $V$.  Note
that a core sphere of $V$ and a core sphere of $U$ are linked
with linking number 1.  

Recall that $n-e = d-1$.  Topologically $V$ has the form $S^{e-1}
\times S^{d-1} \times [0,1]$.  Up to a $C(n)$-bilipschitz
equivalence,
the set $V$ is bilipschitz to a Riemannian product of the
following form: $Core \times E[d-1]
\times [0, E_d]$.  Here $Core$ is a copy of $S^{e - 1}$ equipped
with an ellipsoidal metric with principal axes $E_{n-e+1} \le ...
\le E_n$.

We define the map $F$ on $V$ as follows.  We let $\pi_V$ be the
projection from $V$ to $E[d-1] \times [0, E_d]$.  Next, there is
a map $\psi_V$ from $E[d-1] \times [0, E_d]$ to $E[d]$, taking
the two boundary components of the domain to the two tips of the
range.  The map has degree 1 and Lipschitz constant at most
$C(n)$.  We define $F$ on $V$ to be the composition $\psi_V \circ
\pi_V$.  The map $\pi_V$ takes the boundary of $V$ to the
boundary of $E[d-1] \times [0, E_d]$, and so the map $F$ takes
the boundary of $V$ to the tips of $E[d]$.  Now, the
identification map $E[d] \rightarrow X$ takes the tips of $E[d]$
to the basepoint of $X$.  Therefore, we can think of $F$ as a map
from $V$ to $X$ taking the boundary of $V$ to the basepoint of
$X$.

We have now defined $F$ on $U$ and on $V$.  The sets $U$ and $V$
are disjoint, and $F$ maps their boundaries to the basepoint of
$X$.  We extend $F$ to all of $E$ by mapping the rest of $E$ to
the basepoint of $X$.  The tips of $E$ are in the complement of
$U \cup V$, and so they get mapped to the basepoint of $X$ as
claimed.  The map $F$ has Lipschitz constant at most $C(n)$.

The last step is to check that the linking invariant of $F$ is
equal to 1.  We let $q_1$ be a generic point in $S^d \subset X$
and $q_2$ a generic point in $S^e \subset X$.  The preimage
$F^{-1}(q_1)$ is a core sphere of $V$.  The preimage
$F^{-1}(q_2)$ is a core sphere of $U$.  These two core spheres
have linking number 1.
\endproof

Now we return to the proof of Theorem 2.  First we apply the
lemma with $e = k_1$.  We get a map from $E$ to $X$.  Recall that
$X$ is formed by gluing together $S^{n-k_1+1}$ and $S^{k_1}$. 
The metric on $S^{n-k_1+1}$ is $E[n-k_1+1]$.

The next step of the proof is to apply the lemma again with
domain $E[n-k_1+1]$.  This time, we choose $e=k_2$.  The lemma
gives us a map $F_2$ from $E[n-k_1+1]$ to a space $X'$.  The
space $X'$ is formed from $S^{k_3} \cup S^{k_2}$ by identifying
the basepoint of $S^{k_2}$ and two antipodal points of $S^{k_3}$. 
The map of $F_2$ sends the tips of $E[n-k_1+1]$ to the basepoint
of $X'$.  Therefore, $F_2$ extends to a map from $X$ to $X' \vee
S^{k_1}$, taking the copy of $S^{k_1}$ in $X$ identically to the
copy of $S^{k_1}$ in $X' \vee S^{k_1}$.  By composing $F_2 \circ
F$, we get a map from $E$ to $X' \vee S^{k_1}$.  We call this map
$\Phi_1$, and we abbreviate $Y = X' \vee S^{k_1}$ 

The second order linking invariant is defined for a map $\Phi$
from $S^n$ to $Y$ in the usual way.  Namely, let $q_i$ be a
generic point in $S^{k_i} \subset Y$, and repeat the usual
procedure with the fibers $\Phi^{-1}(q_i)$.  The map $\Phi_1$ has
$L_2(\Phi_1) = 1$.  (This calculation is essentially the same
as the calculation of $L_2(g^+ \circ f)$ from Section 1.)

The space $Y$ is equipped with a metric $g$, and with respect to
this metric the map $\Phi_1$ has Lipschitz constant at most
$C(n)$.  The metric on $S^{k_1} \subset Y$ is the one-point
compactification of the rectangle with dimensions $E_{n-k_1+1}
\times ... \times E_n$.  The metric on $S^{k_2} \subset Y$ is the
one-point compactification of the rectangle with dimensions
$E_{n-k_1-k_2+2} \times ... \times E_{n-k_1+1}$.  The metric on
$S^{k_3} \subset Y$ is the ellipsoidal metric $E[k_3]$.

Next we construct a map $\alpha$ from $Y$ to $S^{k_3} \vee S^{k_2}
\vee S^{k_1}$, which takes $S^{k_i} \subset Y$ to $S^{k_i}$.  We
put the standard unit sphere metric on each sphere in the range. 
For large $L$, we can find $\alpha$ with Lipschitz constant $L$
and with degree $D_1$ at least $c(n) E_{n-k_1+1} ... E_n L^{k_1}$
on $S^{k_1}$, with degree $D_2$ at least $c(n) E_{k_3} ...
E_{n-k_1+1} L^{k_2}$ on $S^{k_2}$, and with degree $D_3$ at least
$c(n) E_1 ... E_{k_3} L^{k_3}$ on $S^{k_3}$.

The map $\Phi$ is $\alpha \circ \Phi_1$.  It has Lipschitz
constant at most $C(n) L$ and $L_2(\Phi) = D_1 D_2 D_3$, which is
at least $c(n) E_{k_3} E_{n-k_1+1} Vol(E) L^{n+2}$. \endproof

\section{Application to k-dilation of degree non-zero maps}

Our two theorems immediately imply a new lower bound on the
k-dilation of a map from one ellipse to another.

\begin{theorem 3} Let $E, E'$ be n-dimensional ellipses. 
Let $E_0 \le ... \le E_n$ be the principal axes of $E$.  Let
$E_0' \le ... \le E_n'$ be the principal axes of $E'$.  Let $Q_i
= E_i'/E_i$.  Suppose that $\Phi$ is a map from $E$ to $E'$ with
degree $D$.  Suppose that $2 \le k_1 \le k_2 \le k_3$, $n+2 = k_1
+ k_2 + k_3$ and $k \le k_1$.  Then the following inequality
holds.

$$Dil_k(\Phi) > c(n) [|D| Q_{n-k_1+1} Q_{k_3} Q_1 ...
Q_n]^{\frac{k}{n+2}}.$$

\end{theorem 3}

\proof By Theorem 2, we can
find a map $F$ from $E'$ to $S^{k_1}
\vee S^{k_2} \vee S^{k_3}$ with Lipschitz constant $L$ large and
$L_2(F) \ge c(n) E_{n-k_1+1}' E_{k_3}' E_1' ... E_n' L^{n+2}$. 
Therefore, $|L_2(F \circ \Phi)|$ is at least $c(n) |D|
E_{n-k_1+1}' E_{k_3}' E_1' ... E_n' L^{n+2}$.  On the other hand,
the map $F \circ \Phi$ has $k_1$-dilation at most
$Dil_k(\Phi)^{k_1/k} L^{k_1}$.  By Theorem 1, the norm of $L_2(F
\circ \Phi)$ must be at most $C(n) E_{n-k_1+1} E_{k_3} E_1 ... E_n
L^{n+2} Dil_k(\Phi)^{(n+2)/k}$.  Comparing the upper and lower
bounds we get the estimate. \endproof

\end{document}